\documentclass[pamm,a4paper,fleqn]{w-art}
\usepackage{times,cite,w-thm}
\usepackage[T1]{fontenc}
\usepackage[utf8]{inputenc}
\theoremstyle{plain}

\theoremstyle{definition}

\usepackage{graphicx}

\usepackage{amssymb}
\usepackage{bm}
\usepackage{tcolorbox}
\newcommand{\mR}{\mathbb{R}}

\newcommand{\SO}{\mathit{SO}}
\newcommand{\inv}{^{-1}}
\newcommand{\vzero}{\mathbf{0}}
\newcommand{\diff}[1][]{\mathrm{d}#1}

\newcommand{\T}{^{\mathop{\mathrm{T}}}}

\newcommand{\diag}{{\mathop{\mathrm{diag}}}}
\newcommand{\blkdiag}{{\mathop{\mathrm{blkdiag}}}}

\newcommand{\Exp}{\mathop{\mathrm{Exp}}}
\newcommand{\dExp}{\mathop{\mathrm{dExp}}}
\newcommand{\dExpInv}{\mathop{\mathrm{dExp}^{+}}}

\newcommand{\vga}{\bm{\gamma}}

\newcommand{\vka}{\bm{\kappa}}
\newcommand{\vla}{\bm{\lambda}}

\newcommand{\vta}{\bm{\tau}}

\newcommand{\vvph}{\bm{\phi}}

\newcommand{\vTh}{\mathbf \Theta}

\newcommand{\vOm}{\mathbf \Omega}

\newcommand{\va}{\mathbf a}
\newcommand{\vb}{\mathbf b}
\newcommand{\vc}{\mathbf c}

\newcommand{\ve}{\mathbf e}
\newcommand{\vf}{\mathbf f}
\newcommand{\vg}{\mathbf g}

\newcommand{\vl}{\mathbf l}
\newcommand{\vm}{\mathbf m}
\newcommand{\vn}{\mathbf n}

\newcommand{\vp}{\mathbf p}
\newcommand{\vq}{\mathbf q}
\newcommand{\vr}{\mathbf r}
\newcommand{\vs}{\mathbf s}

\newcommand{\vu}{\mathbf u}
\newcommand{\vv}{\mathbf v}

\newcommand{\vA}{\mathbf A}
\newcommand{\vB}{\mathbf B}

\newcommand{\vE}{\mathbf E}

\newcommand{\vI}{\mathbf I}

\newcommand{\vK}{\mathbf K}

\newcommand{\vM}{\mathbf M}

\newcommand{\vS}{\mathbf S}

\newcommand{\vW}{\mathbf W}

\usepackage{pgfplots}
\usepgfplotslibrary{groupplots}
\pgfplotsset{compat=1.18}
\usepackage{siunitx}
\definecolor{myred}{RGB}{238,102,119}
\definecolor{mygreen}{RGB}{34,136,51}
\definecolor{myblue}{RGB}{68,119,170}

\begin{document}
%% \def\leftmark{Session title}
%%
%%    The information for the title page will be placed between
%%    \begin{document} and \maketitle. The order of most entries
%%    is determined by the class file and can not be changed by
%%    rearranging them. The maketitle command follows after the
%%    abstract.
%%
%%    Most of the following commands will be completed by the publisher.
%%
%%    \renewcommand{\copyrightyear}{2016}
%%    \DOIsuffix{pamm.20161zzzz}
%%    \Volume{16} 
%%    \Year{2016} 
%%    \pagespan{1}{}
%%
%%    The short title is optional:

\TitleLanguage[EN]
\title[The short title]{A Mixed Discrete Cosserat Rod Formulation}

%% Please do not enter footnotes or \inst{}-notes into the optional
%% argument of the author command. 

%% Please delete not needed author entries.
%% Information for the first author.
\author{\firstname{Tianxiang}  \lastname{Dai}\inst{1,}%
  \footnote{Corresponding author: e-mail \ElectronicMail{dai@inm.uni-stuttgart.de},
            phone +49\,711\,685\,68166,
            fax +49\,711\,685\,66277}}

\address[\inst{1}]{\CountryCode[DE] University of Stuttgart, Institute for Nonlinear Mechanics, Stuttgart, Germany}
%%
%%    Information for the second author
\author{\firstname{Marco} \lastname{Herrmann}\inst{2}}

\address[\inst{2}]{\CountryCode[NL] Eindhoven University of Technology, Department of Mechanical Engineering, Eindhoven, The Netherlands}
%%
%%    Information for the third author
\author{\firstname{Jonas} \lastname{Breuling}\inst{1}}
  
\author{\firstname{Remco I.} \lastname{Leine}\inst{1}}
 	
\author{\firstname{Simon R.} \lastname{Eugster}\inst{2}}
%%
%%    \dedicatory{This is a dedicatory.}
%%
%%    Abstract is required.
\AbstractLanguage[EN]
\begin{abstract}
	In this communication we propose a discrete Cosserat rod formulation in which a slender elastic rod is represented as a chain of rigid bodies (nodes) coupled by compliant elastic forces and moments acting between adjacent node pairs. Discrete dilatation, shear, torsion and curvature strain measures are evaluated from the relative kinematics of each node pair, while the constitutive behavior is expressed in compliance form through independent stress degrees of freedom. We show that the resulting model arises rigorously from a mixed Petrov--Galerkin Cosserat rod finite element formulation (FEM) at linear kinematic interpolation order when the internal virtual work is integrated by the midpoint rule and the external and inertial contributions by the trapezoidal rule. The proposed formulation inherits the robustness and the absence of locking from the underlying mixed FEM while simultaneously exposing a two-node coupling structure that mirrors discrete rod models from the computer graphics community. This is in sharp contrast to the dense coupling of strain-parameterized reduced-order models often used in soft robotic applications. Three numerical examples involving piecewise-varying cross sections, tendon-driven actuation under different spatial discretizations, and coupled longitudinal-torsional dynamics confirm the accuracy, robustness, and convergence behavior of the presented approach.
\end{abstract}
\maketitle                   % Produces the title.

\section{Introduction}
Spatial rod theories are widely adopted to model the mechanics of slender structures in fields ranging from soft robotics and structural mechanics to multibody dynamics and computer graphics. For numerical simulation, a variety of modeling strategies have been proposed, including discrete (geometric) rod models that originate in the computer graphics community~\cite{Bergou2008, Kugelstadt2016, Gazzola2018, Hsu2025} and in the multibody community~\cite{Lang2010, Jung2010}, shooting-based formulations of the Cosserat rod boundary value problems (BVPs) developed for continuum robotics~\cite{Rucker2010, Till2019, Tummers2023}, strain-parameterized reduced-order models~\cite{Renda2018, Boyer2021, Mathew2023}, and Cosserat rod finite element formulations (FEMs)~\cite{Harsch2021, Harsch2023, Herrmann2025} from structural mechanics.

In the computer graphics community, discrete rod models~\cite{Bergou2008, Kugelstadt2016, Gazzola2018, Hsu2025} have been widely employed to simulate large-scale systems such as hair strands, yarns, ropes, and tree branches. These methods focus on visual plausibility, real-time performance, and simulation robustness rather than on high-fidelity quantitative accuracy. To meet the demands of fast simulation, a discrete rod is represented as a sequence of vertices and edges, where vertices carry the position of the centerline, while orientation information is associated with the edges. This representation is straightforward to implement, exposes a regular two-node coupling topology that, combined with local relaxation solvers, has been shown to enable parallel GPU implementations with millions of vertices~\cite{Hsu2025}.

In the robotics community, the shooting method for BVPs has become a standard tool for the statics~\cite{Rucker2010} and dynamics~\cite{Till2019} of Cosserat rods. While the shooting method is conceptually simple and easy to implement, its convergence is sensitive to the initial guess of the unknown boundary state. Consequently, challenging loading scenarios typically require continuation through multiple load increments to ensure convergence~\cite{Tummers2023}. As an alternative, strain-based formulations of Cosserat rods have been developed, such as the piecewise constant strain model~\cite{Renda2018}, the geometric variable-strain approach~\cite{Renda2020, Boyer2021}, and FEM-like local strain parameterizations~\cite{Mathew2024, Mathew2025}. By projecting the infinite-dimensional configuration space of the rod to a finite-dimensional one through global or piecewise-local strain bases, these models achieve computationally efficient discrete representations. However, because the strain coordinates encode the relative deformation between adjacent cross sections, the absolute pose of any cross section depends recursively on all preceding strain variables. As a consequence, the mass and damping matrices of the discrete system are dense and the contributions of distributed external loads are highly nonlinear in the generalized coordinates, which has motivated the use of recursive Newton--Euler-type algorithms to evaluate the equations of motion at acceptable cost~\cite{Renda2018, Boyer2021, Mathew2025}. 

For modeling soft robotic systems, Cosserat rod FEMs remain comparatively under-explored~\cite{Armanini2023}, despite their maturity in the structural mechanics literature~\cite{Harsch2023,Herrmann2025}. A recurring obstacle is that any Cosserat rod FEM requires a careful treatment of large rotations, objective kinematic interpolations, and locking-mitigation strategies, all of which raise the barrier to entry and make integration with existing rigid multibody simulation pipelines non-trivial. To address this gap, we propose a novel discrete Cosserat rod model consisting of a sequence of rigid bodies (subsequently called nodes) coupled by compliant elastic forces and moments acting between adjacent node pairs. Discrete strain measures (dilatation, shear, torsion and curvature) describing rod deformation are computed directly from the relative positions and orientations of each node pair, and the constitutive laws relating these discrete strains to the corresponding internal forces and moments are given in compliance form. The resulting model is straightforward to implement and can be integrated into existing rigid body simulation frameworks with minimal additional infrastructure. Importantly, we show that this discrete multibody representation is not an empirical or heuristic construction but emerges rigorously from the mixed Petrov--Galerkin Cosserat rod FEM recently introduced by Herrmann et al.~\cite{Herrmann2025}, in which the resultant internal forces and moments are interpolated as independent fields via the Hellinger--Reissner principle, combined with the Petrov--Galerkin framework of Harsch et al.~\cite{Harsch2023}, where the internal virtual work is integrated by the midpoint rule and the external and inertial contributions by the trapezoidal rule. Our model therefore inherits the robustness and the absence of locking from the underlying mixed FEM, while simultaneously exposing the same two-node kinematic topology as discrete rod models from the computer graphics community~\cite{Bergou2008, Kugelstadt2016, Gazzola2018, Hsu2025}. Adapting GPU-parallel solver strategies of the kind demonstrated in~\cite{Hsu2025} to the present mixed discrete formulation is left for future work. 

The remainder of this paper is organized as follows. Section~\ref{sec:kinematics} introduces the discrete Cosserat rod kinematics, including the nodal parameterization, the discrete strain measures, and the independent stress degrees of freedom arising from the mixed formulation. Section~\ref{sec:eom} presents the equations of motion governing the dynamics of the discrete Cosserat rod model, together with their reduction to the nonlinear static equilibrium equations. Section~\ref{sec:examples} presents three numerical examples to validate the accuracy, convergence behavior, and applicability of the proposed model in both static and dynamic settings; specifically, a two-segment helix with piecewise-varying cross sections, a tendon-driven continuum manipulator with tapered geometry undergoing large deformations, and a Wilberforce pendulum exhibiting coupled longitudinal--torsional oscillations. Finally, Section~\ref{sec:conclusion} summarizes the main findings and outlines directions for future work.

\section{Discrete Cosserat rod kinematics} \label{sec:kinematics} 
Following~\cite{Antman2005}, the Euclidean 3-space $\mathbb{E}^3$ is introduced as an abstract 3-dimensional real inner product space. Given an arbitrary right-handed orthonormal basis $B = \{\ve_x^B, \, \ve_y^B, \, \ve_z^B\}$, a vector $\va = a_x^B \ve_x^B + a_y^B \ve_y^B + a_z^B \ve_z^B \in \mathbb{E}^3$ can be represented in the $B$-basis by the coordinate triple ${}_B \va = (a_x^B , \, a_y^B , \, a_z^B) \in \mR^3$, which collects the scalar components $a_i^B$, $i \in \{x, y, z\}$. This notation clearly distinguishes between the abstract Euclidean space $\mathbb{E}^3$ and the space of real triples $\mR^3$. The coordinate representations ${{}_{B_0} \va}$ and ${{}_{B_1} \va}$ of the same vector $\va$ with respect to two different bases $B_0$ and $B_1$ are related by the proper orthogonal transformation matrix $\vA_{B_0 B_1} = \big({}_{B_0} \ve_x^{B_1},\, {}_{B_0} \ve_y^{B_1},\, {}_{B_0} \ve_z^{B_1}\big) \in \SO(3) = \{\vA \in \mR^{3 \times 3} \mid \vA\T \vA = \vE, \, \det(\vA) = 1\}$ through ${}_{B_0} \va = \vA_{B_0 B_1} \, {}_{B_1} \va$. The zero matrix in $\mR^{3\times 3}$ and the zero $n$-tuple in $\mR^n$ are both denoted $\vzero$ (or $\vzero_n$ when the dimension must be made explicit), and $\vE \in \mR^{3\times 3}$ denotes the identity matrix. To keep expressions concise, basis subscripts are included only when necessary, and the coordinate space of remaining quantities should be inferred from context. The operator $j \colon \mR^3 \to \mathfrak{so}(3) = \{\vS \in \mR^{3\times 3} \mid \vS\T = -\vS\}$ is the linear bijection induced by the vector product in agreement with $j(\va)\vb = \widetilde{\va}\vb = \va \times \vb$ for all $\va, \vb \in \mR^3$.

\begin{figure}[b]
    \centering
    \includegraphics[width=120mm]{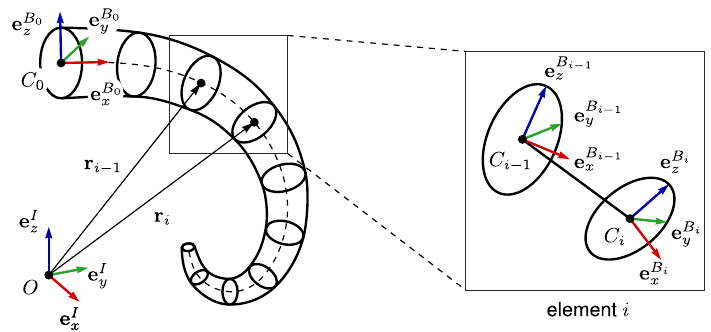}
    \caption{Discrete Cosserat rod kinematics. The rod is discretized into $n$ segments by $n+1$ separating cross sections (nodes) at centerline points $\{C_0, \ldots, C_n\}$. Each node $i$ is parameterized by the rod centerline position $\vr_i$ and has a node-fixed basis $\{\ve_x^{B_i}, \ve_y^{B_i}, \ve_z^{B_i}\}$. The inset highlights two adjacent nodes at $C_{i-1}$ and $C_i$ forming element $i$.}
    \label{fig:discrete_rod}
\end{figure}

Consider a Cosserat rod of reference length $L$, uniform mass density $\rho$, and variable cross-sectional area $A(s)$, where $s \in [0, L]$ denotes the reference arc-length parameter. At time $t$, the configuration of the rod is described by the centerline position $\vr(s, t) := {}_I\vr_{OC}(s, t) \in \mR^3$ and the cross-sectional orientation $\vA(s, t) := \vA_{IB}(s, t) \in \SO(3)$, where $I$ denotes the inertial basis, $O$ the origin, $C(s)$ the centerline point, and $B(s)$ the cross-section-fixed basis. The centerline is assumed to pass through the center of mass of each cross section, with the planar second moments of area $I_y(s)$ and $I_z(s)$, and the polar second moment of area $J_x(s)$, defined about the centroidal axes. For discretization, the rod is divided into $n$ segments, yielding $n_\mathrm{node} = n+1$ separating cross sections located at arc-length positions $\{s_i \mid i = 0, \ldots, n\}$ with $s_0 = 0$ and $s_n = L$. The $i$-th segment $(i \in \{1, \ldots, n\})$ spans interval $[s_{i-1}, s_i]$ and has element length $L_i^\mathrm{el} = s_i - s_{i-1}$. The cross sections are treated as independent rigid bodies (nodes), each having a node-fixed basis $B_i=B(s_i)=\{\ve_x^{B_i}, \ve_y^{B_i}, \ve_z^{B_i}\}$ and an orientation $\vA_i(t) = \vA(s_i, t) \in \SO(3)$. The kinematic structure of the discretization is illustrated in Figure~\ref{fig:discrete_rod}. The cross-sectional area at the $i$-th node is denoted by $A_i = A(s_i)$, while the area at the midpoint of the $i$-th segment is $A_i^\mathrm{el} = A(s_i^\mathrm{el})$, where $s_i^\mathrm{el} = (s_{i-1} + s_i)/2$. Similarly, the polar and planar second moments of area at the node are given by $J_{x,i} = J_x(s_i), I_{y,i} =  I_y(s_i)$ and $I_{z,i} =  I_z(s_i)$, whereas the corresponding element values evaluated at the element midpoint are $J_{x,i}^\mathrm{el} = J_x(s_i^\mathrm{el}), I_{y,i}^\mathrm{el} =  I_y(s_i^\mathrm{el})$ and $I_{z,i}^\mathrm{el} =  I_z(s_i^\mathrm{el})$. The mass of each segment is lumped equally to its two boundary nodes as will be shown to result from the trapezoidal quadrature of the inertial virtual work in Section~\ref{sec:relation_to_FEM}. The effective nodal mass $m_i$ and inertia tensor $\vTh_i\in \mR^{3\times 3}$, expressed with respect to the centerline point $C_i = C(s_i)$, of the $i$-th node are therefore given by
\begin{equation} \label{eq:nodal_mass_inertia}
    \begin{split}
        m_i & = \rho L_i A_i,
        \\
        \vTh_i & =\rho L_i \diag(J_{x,i}, I_{y,i}, I_{z,i}),
    \end{split}
    \qquad \text{with }
    L_i = \begin{cases}
        L_1^\mathrm{el}/2, & i = 0, \\
        (L_i^\mathrm{el} + L_{i+1}^\mathrm{el})/2, & 0 < i < n, \\
        L_n^\mathrm{el}/2, & i = n.
    \end{cases}
\end{equation}
The distributed external force density $\vb(s) \in \mR^3$ (given in the inertial basis $I$) and moment density $\vc(s) \in \mR^3$ (given in the cross-section-fixed basis $B$) are likewise lumped to the nodes as concentrated forces $\vb_i = L_i \vb(s_i)$ and moments $\vc_i = L_i \vc(s_i)$.

Following the Petrov--Galerkin Cosserat rod formulation of~\cite{Harsch2023}, each node is parameterized by the generalized coordinate $\vq_i = (\vr_i, \vp_i) \in \mR^7$, where the quaternion 
$\vp_i = (p_0, \vp_\mathrm{vec}) \in \mR^4$, with scalar part $p_0 \in \mR$ and vector part $\vp_\mathrm{vec} \in \mR^3$, represents the orientation via
\begin{equation} \label{eq:exp_operator}
    \vA_i = \Exp(\vp_i) = \vE + \frac{2}{\|\vp_i\|^2}
    \big(p_0 \,\widetilde{\vp}_\mathrm{vec} + \widetilde{\vp}_\mathrm{vec}^2\big).
\end{equation}
The nodal velocity state is represented using the minimal velocity coordinates $\vu_i = (\vv_i, \vOm_i) \in \mR^6$, where $\vv_i = \dot{\vr}_i$ denotes the centerline velocity, and $\vOm_i$ denotes the angular velocity of the node-fixed basis $B_i$ relative to the inertial basis $I$. The angular velocity is related to the nodal transformation matrix through $\dot{\vA}_i = \vA_i\,\widetilde{\vOm}_i$. Let $\dExp(\vp_i) \in \mR^{3\times 4}$ denote the tangent operator of the quaternion
\begin{equation} \label{eq:tangent_operator}
    \dExp(\vp_i) = \frac{2}{\|\vp_i\|^2}
    \begin{pmatrix}
        -\vp_\mathrm{vec} & p_0 \vE - \widetilde{\vp}_\mathrm{vec}
    \end{pmatrix},
\end{equation}
and let $\dExpInv(\vp_i)\in \mR^{4\times 3}$ be its Moore--Penrose inverse,
\begin{equation}
    \dExpInv(\vp_i) = \frac{1}{2}
    \begin{pmatrix}
        -\vp_\mathrm{vec}\T \\[2pt]
        p_0 \vE + \widetilde{\vp}_\mathrm{vec}
    \end{pmatrix},
\end{equation}
such that $\dExp(\vp_i)\dExpInv(\vp_i) = \vE$, as shown in~\cite{Rucker2018}. The kinematic differential equation relating $\dot{\vq}_i$ to $\vu_i$ can then be written as
\begin{equation}
    \dot{\vq}_i = \begin{pmatrix}
        \dot{\vr}_i \\  \dot{\vp}_i
    \end{pmatrix} = \vB_i(\vq_i)\,\vu_i, 
    \qquad \vB_i(\vq_i) = \begin{pmatrix}
        \vE & \vzero \\ \vzero & \dExpInv(\vp_i)
    \end{pmatrix},
    \label{eq:kin_diff_eq}
\end{equation}
where the rotational kinematics satisfy $\dot{\vp}_i = \dExpInv(\vp_i)\vOm_i$. 
Assembling over all nodes, the global generalized coordinates and velocities are
\begin{equation}
    \vq = (\vq_0, \ldots, \vq_n) \in \mR^{7(n+1)}, \qquad
    \vu = (\vu_0, \ldots, \vu_n) \in \mR^{6(n+1)}.
\end{equation}

The discrete strain measures of the $i$-th segment are defined using the relative kinematics between nodes $(i-1)$ and $i$. The segment orientation is approximated by the mean quaternion $\vp_i^\mathrm{el} = (\vp_i + \vp_{i-1})/2$ with corresponding transformation matrix $\vA_i^\mathrm{el} = \Exp(\vp_i^\mathrm{el})$. The discrete dilatation-shear strain $\vga_i := {}_{B_i^\mathrm{el}} \vga_i \in \mR^3$ and the discrete torsion-curvature strain $\vka_i := {}_{B_i^\mathrm{el}} \vka_i \in \mR^3$, both expressed in the segment basis $B_i^\mathrm{el}$, are defined as
\begin{equation}
    \vga_i = \big(\vA_i^\mathrm{el}\big)\T 
              \frac{\vr_i - \vr_{i-1}}{L_i^\mathrm{el}}, \qquad
    \vka_i = \dExp(\vp_i^\mathrm{el})\,\frac{\vp_i - \vp_{i-1}}{L_i^\mathrm{el}}.
    \label{eq:strain_measures}
\end{equation}
Note that, due to the normalization inherent in the quaternion operators $\Exp$ in~\eqref{eq:exp_operator} and $\dExp$~in \eqref{eq:tangent_operator}, the quantities $\vA_i^\mathrm{el}$, $\vga_i$, and $\vka_i$ remain consistent even if the mean quaternion $\vp_i^\mathrm{el}$ is not of unit length. 
Following the Hellinger--Reissner variational principle underlying the mixed FEM of~\cite{Herrmann2025}, the internal forces $\vn_i := {}_{B_i^\mathrm{el}} \vn_i \in \mR^3$ and moments $\vm_i := {}_{B_i^\mathrm{el}} \vm_i \in \mR^3$ of the $i$-th segment are treated as independent stress-like degrees of freedom. Collecting the stress variables of all segments yields
\begin{equation}~\label{eq:stress}
    \vla_c = (\vn_1, \vm_1, \ldots, \vn_n, \vm_n) \in \mR^{6n}.
\end{equation}
The constitutive law relating the discrete strains~\eqref{eq:strain_measures} to the corresponding internal forces and moments~\eqref{eq:stress} is expressed in compliance form as
\begin{equation}
    \vzero_6 = \vK_{c,i}\inv \begin{pmatrix}
        \vn_i \\ \vm_i
    \end{pmatrix} - \vl_i(\vq),
\end{equation}
where the diagonal stiffness matrix and the deformation measure of the segment are defined by
\begin{align} \label{eq:el_stiff_matrix}
    \vK_{c,i} &= L_i^\mathrm{el}\, 
        \diag\big(EA_i^\mathrm{el},\; GA_i^\mathrm{el},\; GA_i^\mathrm{el},\; 
                     GJ_{x,i}^\mathrm{el},\; EI_{y,i}^\mathrm{el},\; EI_{z,i}^\mathrm{el}\big) \in \mR^{6\times 6}, \\
    \vl_i(\vq) &= L_i^\mathrm{el}\,
        \big(\vga_i - \vga_i^0,\; \vka_i - \vka_i^0\big) \in \mR^6,
\end{align}
where $E$ and $G$ denote the Young's and shear moduli, respectively, and the second moments of area $J_{x,i}^\mathrm{el}$, $I_{y,i}^\mathrm{el}$ and $I_{z,i}^\mathrm{el}$ are defined in Section~\ref{sec:kinematics}. The reference strains $\vga_i^0$ and $\vka_i^0$ are computed from a stress-free reference configuration.
\section{Equations of motion} \label{sec:eom} 
The dynamics of the discrete rod system described above is governed by the following equations of motion
\begin{align}
    \dot{\vq} &= \vB(\vq)\,\vu, \label{eq:kin_diff_eq_sys} \\
    \vM\,\dot{\vu} + \vf^\mathrm{gyr}(\vu) &= \vf^\mathrm{ext}(\vq) 
        + \vW_c(\vq)\,\vla_c, \label{eq:dyn_equilibrium} \\
    \vzero_{6n} &= \vK_c^{-1}\,\vla_c - \vl_c(\vq). \label{eq:compliance_sys}
\end{align}
In the kinematic differential equation~\eqref{eq:kin_diff_eq_sys}, the operator $\vB(\vq)$ applies the operator $\vB_i$ in~\eqref{eq:kin_diff_eq} to each node, such that
\begin{equation}
    \vB(\vq) = \blkdiag \big(\vB_0(\vq_0),\; \ldots, \;\vB_n(\vq_n)\big).
\end{equation}
Together with consistent initialization $\|\vp_i(0)\| = 1$ for all $i$ and normalization after each integration step, the unit length of the quaternions is implicitly maintained. In the kinetic part of the equations of motion~\eqref{eq:dyn_equilibrium}, the system mass matrix
\begin{equation}
    \vM = \blkdiag(m_0 \vE,\; \vTh_0,\; \ldots,\; m_n \vE,\; \vTh_n) \in \mR^{6(n+1)\times 6(n+1)}
\end{equation}
is constant and diagonal, where the nodal mass $m_i$ and inertia tensor $\vTh_i$ are introduced in~\eqref{eq:nodal_mass_inertia}. The gyroscopic force vector is accordingly given by
\begin{equation}
    \vf^\mathrm{gyr}(\vu) = 
    \big(\vzero_3,\; \widetilde{\vOm}_0\,\vTh_0\,\vOm_0,\; \ldots,\; 
           \vzero_3,\; \widetilde{\vOm}_n\,\vTh_n\,\vOm_n\big) \in \mR^{6(n+1)}.
\end{equation}
The generalized external force vector obtained by lumping the distributed loads to the nodes is
\begin{equation}
    \vf^\mathrm{ext}(\vq) = 
    \big(\vb_0,\; \vc_0,\; \vb_1,\; \ldots,\;  \vc_{n-1},\; \vb_n,\; \vc_n\big) + \big(\vf_0,\; \vta_0,\; \vzero, \; \ldots,\; \vzero, \; \vf_n,\; \vta_n\big) 
    \in \mR^{6(n+1)},
    \label{eq:discrete_ext_work}
\end{equation}
where the lumped forces $\vb_i$ and moments $\vc_i$ are defined in Section~\ref{sec:kinematics}, and the additional point forces $\vf_0, \vf_n$ and moments $\vta_0, \vta_n$ at the rod boundaries are incorporated straightforwardly into the lumped quantities. The dependence on $\vq$ accounts for load cases in which $\vb_i$, $\vf_0$, or $\vf_n$ are not constant with respect to the inertial basis $I$, or in which $\vc_i$, $\vta_0$, or $\vta_n$ are not constant with respect to the cross-section-fixed basis $B$, see~\cite{Herrmann2025}. The generalized force direction Jacobian $\vW_c(\vq) \in \mR^{6(n+1)\times 6n}$ encodes the action of each segment's internal forces $\vn_i$ and moments $\vm_i$ on the rod nodes. Because each segment only couples its two adjacent nodes, the force Jacobian $\vW_c$ has the banded block structure
\begin{equation} \label{eq:Wc}
    \vW_c = \begin{pmatrix}
        \vW_{c,1}^{l}  & \vzero & \cdots & \vzero \\
        \vW_{c,1}^{r}  & \vW_{c,2}^{l} & \ddots & \vdots \\
        \vzero & \vW_{c,2}^{r} & \ddots & \vzero \\
        \vdots & \ddots & \ddots & \vW_{c,n}^{l} \\
        \vzero & \cdots & \vzero & \vW_{c,n}^{r}
    \end{pmatrix},
\end{equation}
with the $6\times 6$ sub-blocks
\begin{equation}
    \vW_{c,i}^{l} = \begin{pmatrix}
        \vA_i^\mathrm{el} & \vzero \\[2pt]
        \dfrac{L_i^\mathrm{el}}{2}\,\widetilde{\vga}_i & 
        \vE + \dfrac{L_i^\mathrm{el}}{2}\,\widetilde{\vka}_i
    \end{pmatrix}, \qquad
    \vW_{c,i}^{r} = \begin{pmatrix}
        -\vA_i^\mathrm{el} & \vzero \\[2pt]
        \dfrac{L_i^\mathrm{el}}{2}\,\widetilde{\vga}_i & 
        -\vE + \dfrac{L_i^\mathrm{el}}{2}\,\widetilde{\vka}_i
    \end{pmatrix}.
    \label{eq:Wc_blocks}
\end{equation}
The sub-block $\vW_{c,i}^{l}$ (resp. $\vW_{c,i}^{r}$) describes the contribution of the $i$-th segment's internal forces $\vn_i$ and moments $\vm_i$ to the generalized forces acting on the left adjacent node $i-1$ (resp. the right adjacent node $i$). The detailed derivation of these sub-blocks from the internal virtual work of the mixed Cosserat rod FEM~\cite{Herrmann2025} is provided in Section~\ref{sec:relation_to_FEM}. This banded topology is identical to the two-node connectivity of discrete rod models~\cite{Bergou2008, Jung2010, Kugelstadt2016, Hsu2025}, reflecting the fact that both classes of models share the same kinematic graph. Finally, the system compliance equation~\eqref{eq:compliance_sys} assembles the stiffness matrices and deformation measure vectors in~\eqref{eq:el_stiff_matrix} globally for all elements as
\begin{equation} \label{eq:Kc_lq}
    \vK_c = \blkdiag(\vK_{c,1},\; \ldots, \;\vK_{c,n}), \qquad
    \vl_c(\vq) = \big(\vl_1(\vq),\; \ldots, \;\vl_n(\vq)\big).
\end{equation}

The equations of motion~\eqref{eq:kin_diff_eq_sys}--\eqref{eq:compliance_sys} constitute an index-1 differential-algebraic equation (DAE) system, in which the nonlinearity of the deformation measures $\vl_c(\vq)$ is isolated in the algebraic equation~\eqref{eq:compliance_sys}. In principle, the latter could be solved for the compliance forces $\vla_c$ and substituted directly into the kinetic equation~\eqref{eq:dyn_equilibrium}, yielding a reduced system of ordinary differential equations (ODEs). However, retaining the split form is numerically advantageous, as this operator-splitting strategy improves both robustness and computational efficiency compared to the condensed ODE form, in which the nonlinearity of $\vl_c(\vq)$ would otherwise enter the kinetic equation~\eqref{eq:dyn_equilibrium} directly.

For static problems, the inertial and gyroscopic terms vanish identically. Setting $\vu = \dot{\vu} = \vzero$,  the system reduces to the nonlinear equilibrium equations
\begin{align}
    \vzero_{6(n+1)} &= \vf^\mathrm{ext}(\vq) + \vW_c(\vq)\,\vla_c, 
        \label{eq:static_dyn} \\
    \vzero_{6n} &= \vK_c^{-1}\,\vla_c - \vl_c(\vq), 
        \label{eq:static_comp} \\
    \vzero_{n+1} &= \vg_s(\vq),
        \label{eq:static_quat}
\end{align}
following the mixed Cosserat rod FEM~\cite{Herrmann2025} (equation~(26) therein). The additional constraint~\eqref{eq:static_quat}, which collects the unit-quaternion conditions
\begin{equation}
    \vg_s(\vq) = \big(g_{s,0}(\vq_0),\; \ldots,\; g_{s,n}(\vq_n)\big), 
    \qquad g_{s,i}(\vq_i) = \|\vp_i\|^2 - 1 = 0,
\end{equation}
must be imposed explicitly, since the kinematic differential equation~\eqref{eq:kin_diff_eq_sys} that implicitly maintains them in the dynamic case is no longer active.

\section{Relation to FEM} \label{sec:relation_to_FEM}
We now show that the discrete rod model introduced in Sections~\ref{sec:kinematics} and~\ref{sec:eom} follows directly from the mixed Petrov--Galerkin Cosserat rod FEM of Herrmann et al.~\cite{Herrmann2025} at polynomial degree $p=1$, when the element integrals are evaluated by the \emph{midpoint rule} for the internal virtual work and the \emph{trapezoidal rule} for the external virtual work. The inertial terms are obtained by applying the \emph{trapezoidal rule} to the inertial virtual work derived in the same Petrov--Galerkin framework of~\cite{Harsch2023}.

In the mixed FEM of~\cite{Herrmann2025} with $p = 1$, the rod is divided into $n$ elements. The global virtual displacement vector is $\delta \vs = (\delta \vs_0, \ldots, \delta \vs_n) \in \mR^{6(n+1)}$ with $\delta \vs_i = (\delta \vr_i, \delta \vvph_i)$, where $\delta \vr_i := {}_I \delta \vr_{C_i} \in \mR^3$ and $\delta \vvph_i := {}_B \delta \vvph_{IB_i} \in \mR^3$ are the virtual linear and angular displacements of the $i$-th node, respectively. Within the $i$-th element $[s_{i-1}, s_i]$, the centerline position, quaternion, and virtual displacements are interpolated linearly between the nodal values at $s_{i-1}$ and $s_i$ as
\begin{equation}
	\begin{aligned}
		\vr_i^\mathrm{lin}(s) & = (1-\alpha_i(s)) \vr_{i-1} + \alpha_i(s) \vr_i, & 
		\vp_i^\mathrm{lin}(s)  &= (1-\alpha_i(s)) \vp_{i-1} + \alpha_i(s) \vp_i, \\
		\delta \vr_i^\mathrm{lin}(s) & = (1-\alpha_i(s)) \delta \vr_{i-1} + \alpha_i(s) \delta \vr_i, \qquad&
        \delta \vvph_i^\mathrm{lin}(s)  & = (1-\alpha_i(s)) \delta \vvph_{i-1} + \alpha_i(s) \delta \vvph_i,
	\end{aligned}
\end{equation}
where $\alpha_i(s) = (s - s_{i-1})/L_i^\mathrm{el}$. The resultant internal force $\vn_i$ and moment $\vm_i$ are constant within each element (degree $p - 1 = 0$). The strain variables of the $i$-th element are computed with the linear interpolations as
\begin{equation} \label{eq:lin_interp_strain_values}
    \begin{aligned}
        \vga_i^\mathrm{lin}(s) & = \big(\Exp \big(\vp_i^\mathrm{lin}\big)\big   )\T \big(\vr_i^\mathrm{lin} \big)_{,s}  , \qquad
        \vka_i^\mathrm{lin}(s) = \dExp \big(\vp_i^\mathrm{lin}(s)\big) \big(\vp_i^\mathrm{lin} \big)_{,s},
    \end{aligned}
\end{equation}
with $\vga_i^0(s)$ and $\vka_i^0(s)$ denoting the reference strains computed from a stress-free reference configuration.

\paragraph{Internal forces and moments}
The internal virtual work in the mixed Cosserat rod FEM~\cite{Herrmann2025} reads, element-wise,
\begin{equation} \label{eq:cont_int_virt_work}
\begin{split} 
    \delta W^\mathrm{int} = \sum_{i=1}^{n} \int_{s_{i-1}}^{s_i} \bigg\{
    &-\big(\delta \vr_i^\mathrm{lin}\big)_{,s}\T \vA_i^\mathrm{lin}\, \vn_i
    - \big(\delta \vvph_i^\mathrm{lin}\big)_{,s}\T \vm_i
    + \big(\delta \vvph_i^\mathrm{lin}\big)\T \big(\vga_i^\mathrm{lin} \times \vn_i + \vka_i^\mathrm{lin} \times \vm_i\big) \\
    &+ \delta \vn_i\T \big(\vK_{\vga}^{-1}(s)\,\vn_i - (\vga_i^\mathrm{lin} - \vga_i^0)\big)
    + \delta \vm_i\T \big(\vK_{\vka}^{-1}(s)\,\vm_i - (\vka_i^\mathrm{lin} - \vka_i^0)\big)
    \bigg\} \diff s,
\end{split}
\end{equation}
where $\vA_i^\mathrm{lin}(s) = \Exp(\vp_i^\mathrm{lin}(s))$. The dilatation-shear stiffness matrix is given by $\vK_{\vga}(s) = \diag(EA(s), GA(s), GA(s))$ and the torsion-curvature stiffness matrix by $\vK_{\vka}(s) = \diag(GJ_x(s), EI_y(s), EI_z(s))$, where $E$ and $G$ denote the Young's and shear moduli, respectively. Furthermore, $I_y (s)$ and $I_z (s)$ are the planar second moments of area, and $J_x(s)$ is the polar second moment of area, as defined in Section~\ref{sec:kinematics}.

Applying the \emph{midpoint rule} (equivalent to the one-point Gaussian quadrature rule) to evaluate each element integral at the midpoint $s_i^\mathrm{el} = (s_{i-1} + s_i)/2$,
\begin{equation}
    \int_{s_{i-1}}^{s_i} f_i(s)\,\diff s \;\approx\; f_i(s_i^\mathrm{el})\, L_i^\mathrm{el},
\end{equation}
the linear interpolation of the orientation field and its corresponding virtual angular displacement evaluate at the midpoint to
\begin{equation} \label{eq:mid_value_quat}
    \vp_i^\mathrm{lin}(s_i^\mathrm{el}) = \frac{1}{2}(\vp_{i-1} + \vp_i) = \vp_i^\mathrm{el}, \qquad
    \vA_i^\mathrm{lin}(s_i^\mathrm{el}) = \Exp(\vp_i^\mathrm{el}) = \vA_i^\mathrm{el}, \qquad
    \delta \vvph_i^\mathrm{lin}(s_i^\mathrm{el}) = \frac{1}{2}(\delta \vvph_{i-1} + \delta \vvph_i),
\end{equation}
and the arc-length derivatives of the linearly interpolated fields are constant within each element
\begin{equation} \label{eq:mid_value_virt_displacements_derivative}
    \begin{aligned}
        \big(\vr_i^\mathrm{lin}\big)_{,s} &= \frac{\vr_i - \vr_{i-1}}{L_i^\mathrm{el}}, &
        \big(\vp_i^\mathrm{lin}\big)_{,s} &= \frac{\vp_i - \vp_{i-1}}{L_i^\mathrm{el}}, \\
        \big(\delta \vr_i^\mathrm{lin}\big)_{,s} &= \frac{\delta \vr_i - \delta \vr_{i-1}}{L_i^\mathrm{el}}, \qquad&
        \big(\delta \vvph_i^\mathrm{lin}\big)_{,s} &= \frac{\delta \vvph_i - \delta \vvph_{i-1}}{L_i^\mathrm{el}}.
    \end{aligned}
\end{equation}
Evaluating \eqref{eq:lin_interp_strain_values} at the midpoint gives
\begin{equation} \label{eq:mid_value_strain}
    \vga_i^\mathrm{lin}(s_i^\mathrm{el}) = (\vA_i^\mathrm{el})\T \frac{\vr_i - \vr_{i-1}}{L_i^\mathrm{el}} = \vga_i, \qquad
    \vka_i^\mathrm{lin}(s_i^\mathrm{el}) = \dExp(\vp_i^\mathrm{el})\,\frac{\vp_i - \vp_{i-1}}{L_i^\mathrm{el}} = \vka_i,
\end{equation}
which coincide exactly with the discrete strain measures defined in~\eqref{eq:strain_measures}. Substituting all midpoint evaluations~\eqref{eq:mid_value_quat}, \eqref{eq:mid_value_virt_displacements_derivative}, and \eqref{eq:mid_value_strain} into the internal virtual work~\eqref{eq:cont_int_virt_work} and collecting terms by node, one obtains
\begin{equation}
\begin{split}
    \delta W^\mathrm{int} & \approx \sum_{i=1}^{n}
    \bigg\{
    \begin{pmatrix} \delta \vr_{i-1} \\ \delta \vvph_{i-1} \end{pmatrix}\T
    \vW_{c,i}^{l} (\vq)
    \begin{pmatrix} \vn_i \\ \vm_i \end{pmatrix}
    +
    \begin{pmatrix} \delta \vr_i \\ \delta \vvph_i \end{pmatrix}\T
    \vW_{c,i}^{r} (\vq)
    \begin{pmatrix} \vn_i \\ \vm_i \end{pmatrix} 
	% \\ & \phantom{{}\approx \sum_{i=1}^{n} \bigg\{}
    + \begin{pmatrix} \delta \vn_i \\ \delta \vm_i \end{pmatrix}\T
    \left[ \vK_{c,i}^{-1}
    \begin{pmatrix} \vn_i \\ \vm_i \end{pmatrix}
    - \vl_i(\vq) \right]
    \bigg\}, \\
    & = \delta \vs \T \vW_c(\vq) \vla_c + \delta \vla_c \T \big(\vK_c\inv \vla_c - \vl_c(\vq)\big),
\end{split}
\end{equation}
where the generalized force direction Jacobian $\vW_c(\vq)$, the global compliance matrix $\vK_c\inv$, and the global deformation measure $\vl_c(\vq)$ are assembled from the element contributions as given in~\eqref{eq:Wc_blocks} and~\eqref{eq:Kc_lq}.

\paragraph{External forces and moments}
The external virtual work is
\begin{equation}
    \delta W^\mathrm{ext} = \sum_{i=1}^{n} \int_{s_{i-1}}^{s_i} \Big\{
    (\delta \vr_i^\mathrm{lin})\T \vb(s) + (\delta \vvph_i^\mathrm{lin})\T \vc(s)
    \Big\} \diff s + \delta \vr_0\T \vf_0 + \delta \vvph_0\T \vta_0 
    + \delta \vr_n\T \vf_n + \delta \vvph_n\T \vta_n,
\end{equation}
where $\vb(s) \in \mR^3$ is the distributed force density in the inertial basis and $\vc(s) \in \mR^3$ is the distributed moment density in the cross-section-fixed basis, both given per unit reference arc-length, and $\vf_0, \vf_n$ and $\vta_0, \vta_n$ are the point forces and moments applied at the rod boundaries $s=0$ and $s=L$, respectively~\cite{Harsch2023, Herrmann2025}. Applying the \emph{trapezoidal rule} to each element integral,
\begin{equation}
    \int_{s_{i-1}}^{s_i} f_i(s)\,\diff s \;\approx\; \frac{L_i^\mathrm{el}}{2}\Big\{f_i(s_{i-1}) + f_i(s_i)\Big\},
\end{equation}
and using the fact that the linearly interpolated virtual displacement fields evaluate to the nodal values at the endpoints, i.e., $\delta \vr_i^\mathrm{lin}(s_{i-1}) = \delta \vr_{i-1} $ and $\delta \vr_i^\mathrm{lin}(s_{i}) = \delta \vr_{i} $, the \emph{trapezoidal rule} gives
\begin{equation}
    \begin{split}
        \delta W^\mathrm{ext} \approx  \sum_{i=1}^{n} & \bigg\{
            \delta \vr_{i-1}\T \vb(s_{i-1}) + \delta \vvph_{i-1}\T \vc(s_{i-1})
            + \delta \vr_i\T \vb(s_i) + \delta \vvph_i\T \vc(s_i)
            \bigg\} \frac{L_i^\mathrm{el}}{2} \\
            & + \delta \vr_0\T \vf_0 + \delta \vvph_0\T \vta_0 
            + \delta \vr_n\T \vf_n + \delta \vvph_n\T \vta_n.
        \end{split}
\end{equation} 
Rearranging the index summation by collecting the contributions associated with each node $i \in \{0, \ldots, n\}$, and adopting the boundary conventions $L_0^\mathrm{el} = L_{n+1}^\mathrm{el} := 0$ such that the end nodes receive only the contribution from their adjacent element, yields 
\begin{equation}
    \begin{split}
        \delta W^\mathrm{ext} & = \delta \vs \T \vf^\mathrm{ext}(\vq) = \sum_{i=0}^{n} \Big\{
            \delta \vr_i\T \vb(s_i) + \delta \vvph_i\T \vc(s_i)
            \Big\} \frac{L_{i+1}^\mathrm{el} + L_i^\mathrm{el}}{2} + \delta \vr_0\T \vf_0 + \delta \vvph_0\T \vta_0 
            + \delta \vr_n\T \vf_n + \delta \vvph_n\T \vta_n,
    \end{split}
\end{equation}
where the corresponding generalized external force vector results in 
\begin{equation}
    \vf^\mathrm{ext}(\vq) = \big(L_0\vb(s_0),\; L_0\vc(s_0),\; \ldots,\; L_n\vb(s_n),\; L_n\vc(s_n)\big) + \big(\vf_0,\; \vta_0,\; \vzero, \; \ldots,\; \vzero, \; \vf_n,\; \vta_n\big) ,
\end{equation}
which is in agreement with~\eqref{eq:discrete_ext_work}.

\paragraph{Inertial and gyroscopic terms}
The inertial virtual work of the Cosserat rod in the Petrov--Galerkin framework of~\cite{Harsch2023} (equation~(44) therein) reads
\begin{equation}
    \delta W^\mathrm{dyn} = -\sum_{i=1}^{n} \int_{s_{i-1}}^{s_i}
    \bigg\{
    (\delta \vr_i^\mathrm{lin})\T A_{\rho}(s)\,\dot{\vv}(s)
    + (\delta \vvph_i^\mathrm{lin})\T \big(\vI_{\rho}(s)\,\dot{\vOm}(s) + \widetilde{\vOm}(s)\,\vI_{\rho}(s)\,\vOm(s)\big)
    \bigg\} \diff s,
\end{equation}
where $A_{\rho}(s)$ is the mass per unit reference length and $\vI_\rho(s)$ is the cross-sectional inertia tensor, both evaluated at arc-length $s$. Since the centerline passes through the cross-sectional center of mass, the first moment $\vS_{\rho}$ vanishes and the coupling between translational and rotational inertia disappears~\cite{Harsch2023}. Applying the \emph{trapezoidal rule} and rearranging the index summation as before yields
\begin{equation}
    \delta W^\mathrm{dyn} \approx -\sum_{i=0}^{n} \bigg\{
    \delta \vr_i\T A_{\rho}(s_i)\,\dot{\vv}_i
    + \delta \vvph_i\T \big(\vI_{\rho}(s_i)\,\dot{\vOm}_i + \widetilde{\vOm}_i\,\vI_{\rho}(s_i)\,\vOm_i\big)
    \bigg\} L_i,
\end{equation}
which identifies the lumped nodal mass $m_i = \rho L_i A_i = A_{\rho}(s_i)\,L_i$ and nodal inertia tensor $\vTh_i = \vI_\rho(s_i)\,L_i$. In compact form,
\begin{equation}
    \delta W^\mathrm{dyn} = -\delta \vs\T \big(\vM\,\dot{\vu} + \vf^\mathrm{gyr}(\vu)\big),
\end{equation}
with the diagonal mass matrix $\vM = \blkdiag(m_0\vE,\,\vTh_0,\,\ldots,\,m_n\vE,\,\vTh_n)$ and gyroscopic force vector $\vf^\mathrm{gyr}$ as given in Section~\ref{sec:eom}, in agreement with~\eqref{eq:dyn_equilibrium}.

\begin{tcolorbox}[colback=gray!15, colframe=gray!55, boxrule=1pt, arc=5pt]
\paragraph{Summary}
Applying the principle of virtual work $\delta W^\mathrm{int} + \delta W^\mathrm{ext} = 0$ to the midpoint- and trapezoidal-integrated contributions, and requiring the result to hold for all admissible nodal virtual displacements $\delta \vs$ and stress variations $\delta \vla_c$, one recovers precisely the nonlinear static equilibrium~\eqref{eq:static_dyn}--\eqref{eq:static_comp}. The dynamic extension~\eqref{eq:dyn_equilibrium} follows by additionally including the trapezoidal-integrated inertial virtual work, i.e., $\delta W^\mathrm{int} + \delta W^\mathrm{ext} + \delta W^\mathrm{dyn} = 0$. This confirms that the proposed discrete multibody model is not an empirical or heuristic construction but emerges rigorously from existing mixed Petrov--Galerkin Cosserat rod FEMs~\cite{Harsch2023, Herrmann2025}.
\end{tcolorbox}
\section{Numerical examples} \label{sec:examples}
All numerical examples and the simulation code are available in the GitHub repository~\cite{GitRepoDiscreteRodPAMM2026}.
\paragraph{Two helical segments}
To demonstrate the capability of the proposed discrete rod model to handle variable cross sections, we consider an extension of the helix benchmark of~\cite{Herrmann2025}. Rather than having a uniform cross section, the rod in this example consists of two segments with different circular cross-sectional radii, leading to a helical shape for each segment in the deformed configuration. The first segment spans $s \in [0,\,L_1]$ with $L_1=2L/3$ and circular cross-sectional radius $r_1 = L/(2\sigma)$, where $\sigma = 10^2$ is the slenderness ratio. The second segment spans $s \in (L_1,\,L]$ with radius $r_2 = {r_1}/{2^{1/4}}$, giving the cross-sectional area ratio $\frac{A_1}{A_2} = \sqrt{2}$. As a result, the two helices have different curvatures. The target configuration of the first helix has $n_c = 2$ coils along the $\ve_z^I$ axis, with height $h = 40$ and coil radius $R_1 = 10$. The second target helix has the same number of coils $n_c$, but with height $h/2$ and coil radius $R_2 = R_1/2$. The target helical centerline position $\vr^*(s) = {}_I\vr^*_{OC}(s)$ is parametrized by
\begin{equation}
    \begin{split}
        \vr^*(s) & = 
        \begin{cases}
            R_1\big(-\sin\alpha(s),\,\cos\alpha(s),\,c\,\alpha(s)\big), & s\in[0, L_1],\\[2pt]
            R_2\big(-\sin\beta(s),\,\cos\beta(s),\,c\,\beta(s)\big) + \big(0,\, R_2,\, h\big) , & s\in(L_1, L],
        \end{cases}\\[2pt]
        \alpha(s) & = \frac{3\pi n_c}{L}\,s, \qquad \beta(s) = \frac{6\pi n_c}{L}\,(s-L_1), 
    \end{split}
\end{equation}
where $c = h/(2\pi R_1 n_c)$ is the pitch-to-perimeter ratio and the total arc-length is $L = 3\pi R_1 n_c\sqrt{1+c^2}$. The Young's and shear moduli are $E = 1$ and $G = 0.5$, respectively.

The rod is initially straight, clamped at $s = 0$ and loaded at the free end $s = L$ by the terminal moment ${}_B\vta = (\tau_x^B, \; \tau_y^B, \; \tau_z^B)=  \big(c\,k_{t}, \;  0, \; k_{b}\big) / (R_1(1+c^2))$, where $k_{t} = 2GI_1$ and $k_{b} = EI_1$ are the torsional and bending stiffnesses of the first helical segment. The static equilibrium problem~\eqref{eq:static_dyn}--\eqref{eq:static_quat} is solved using a Newton--Raphson scheme with load incrementation, in which the terminal moment is ramped linearly from zero to its prescribed value over 10 equal load steps. Convergence is declared when the absolute and relative residual norms fall below an absolute tolerance of $10^{-10}$ and a relative tolerance of $10^{-6}$, respectively.
\begin{figure}[htbp]
% ================= LEFT IMAGE =================
\begin{minipage}[c]{0.38\textwidth}
    \includegraphics[width=\textwidth]{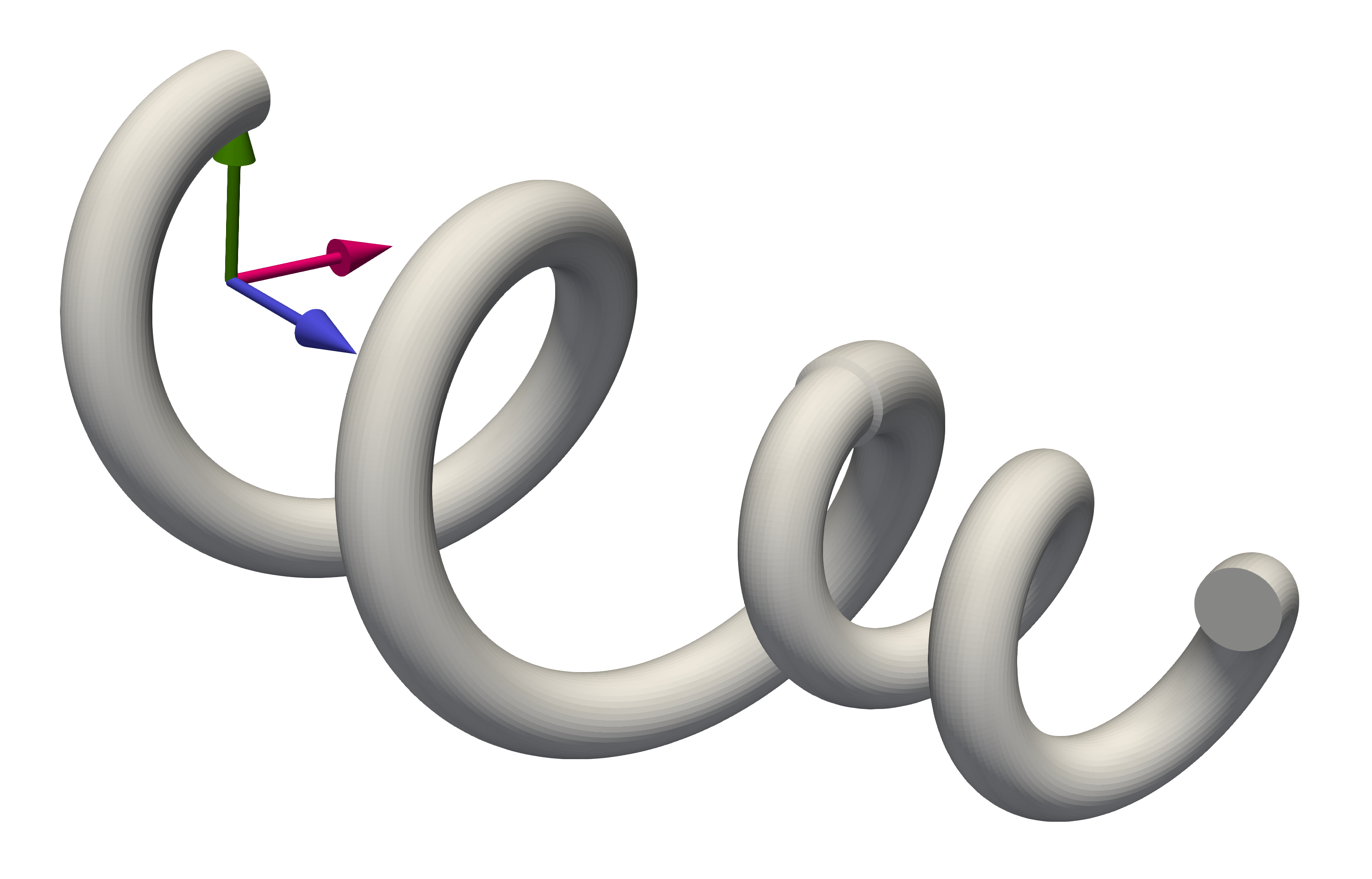}
\end{minipage}
\hfill
% ================= RIGHT PLOT =================
\begin{minipage}[c]{0.61\textwidth}
    \centering
    \includegraphics{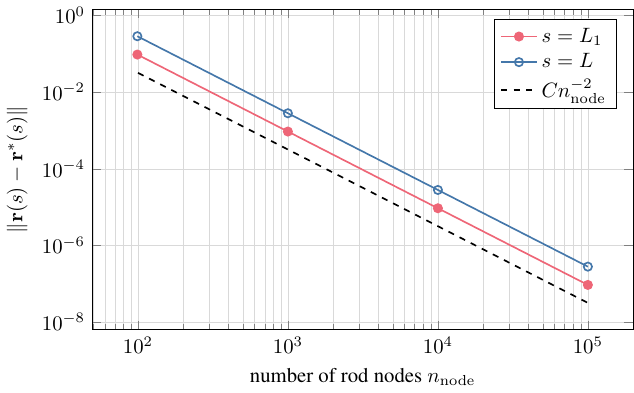}
%    \begin{tikzpicture}
%    \begin{axis}[
%        width=\textwidth,
%        height=7cm,
%        grid=both,
%        grid style={line width=0.2pt, draw=gray!30},
%        tick align=inside,
%        xlabel={number of rod nodes $n_\mathrm{node}$},
%        ylabel={$\lVert \vr(s) - \vr^*(s)\rVert$},
%        xmode=log,
%        ymode=log,
%        % enlarge x limits={abs=10},
%        % enlarge y limits={abs=0.1},
%        legend pos=north east,
%        legend cell align=left
%    ]
%
%    \addplot[
%        myred,
%        thick,
%        mark=*
%    ]
%    table[
%        col sep=comma,
%        x=nnodes,
%        y=err_pos_mid
%    ] {figures/data_double_helix.csv};
%    \addlegendentry{$s=L_1$}
%
%    \addplot[
%        myblue,
%        thick,
%        mark=o
%    ]
%    table[
%        col sep=comma,
%        x=nnodes,
%        y=err_pos_tip
%    ] {figures/data_double_helix.csv};
%    \addlegendentry{$s=L$}
%
%    % dashed ref line
%    \addplot[
%        black,
%        dashed,
%        thick,
%        domain=1e2:1e5
%    ]
%    {316 * x^(-2)}; % 0.095,0.2855
%    \addlegendentry{$Cn_\mathrm{node}^{-2}$}
%    \end{axis}    
%    \end{tikzpicture}
\end{minipage}
\caption{Helix with two different cross-sectional areas. Left: deformed configuration of the rod consisting of two segments with cross-sectional radii $r_1$ and $r_2 = r_1/2^{1/4}$. The reduced bending stiffness of the second segment leads to a higher curvature than in the first segment. Right: convergence of the centerline position error $\lVert \vr(s) - \vr^*(s)\rVert$ at the segment junction $s=L_1$ and at the free end $s=L$ as a function of the number of nodes. The dashed line labeled $Cn_\mathrm{node}^{-2}$ indicates the reference slope corresponding to second-order convergence.} 
\label{fig:helix}
\end{figure}
The deformed configuration and the convergence of the helical centerline position error are shown in Figure~\ref{fig:helix}. Since the target centerline $\vr^*(s)$ is known analytically, the accuracy of the discrete solution can be assessed directly. Both measurement points, at the segment junction $s = L_1$ and at the free end $s=L$, exhibit consistent convergence as the number of nodes increases, confirming the expected approximation properties of the discrete rod formulation. The effect of the reduced bending stiffness in the second segment is instead reflected geometrically. Since $EI_2 = EI_1/2$, the curvature of the second segment is twice that of the first segment under the same moment. Consequently, although the second segment occupies only one third of the total arc-length, it winds two additional coils, so that the deformed rod exhibits a total of 4 coils rather than the 3 coils that a uniform rod of stiffness $EI_1$ would form, as clearly visible in Figure~\ref{fig:helix}. This result confirms that the discrete rod formulation handles piecewise-varying cross sections correctly, with the element-wise compliance law~\eqref{eq:compliance_sys} naturally encoding the local stiffness at each element without any special treatment at the cross-sectional transition $s = L_1$.

\paragraph{Tendon-driven continuum manipulator}
To demonstrate the capability of the proposed discrete rod model for tendon-driven continuum manipulators, we consider a tapered elastic rod actuated by a single routed tendon. The rod has reference length $L = \SI{200}{\milli\meter}$ and a linearly varying circular cross-sectional radius
\begin{equation}
    r(s) = r_0\!\left(1 - \frac{s}{L}(1 - \eta)\right),
    \qquad \eta = \frac{r_\mathrm{tip}}{r_0} = 0.4,
\end{equation}
with base radius $r_0 = \SI{10}{\milli\meter}$ and tip radius $r_\mathrm{tip} = \eta\,r_0 = \SI{4}{\milli\meter}$. The material parameters are the Young's modulus $E = \SI{7e5}{\pascal}$ and shear modulus $G=\SI{2e5}{\pascal}$. The rod is clamped at its bottom end $s = 0$, initially straight and aligned along the $\ve_z^I$ axis. To assess the convergence of the proposed discrete formulation, the equilibrium problem defined by~\eqref{eq:static_dyn}--\eqref{eq:static_quat} is solved using four different spatial discretizations with $n=10$, $50$, $250$, and $1250$ rod elements. The tendon is routed along the lateral surface of the rod. At each node $i$, the tendon attachment point is offset from the centerline by the offset vector $r(s_i)\,\ve_y^{B_i}$ in the body frame. A Newton--Raphson method with load incrementation is employed, where the tendon tension is increased linearly from zero to the maximum value $\lambda_t = \SI{4}{\newton}$ over eight equal load steps. The convergence tolerances are identical to those of the first example, with an absolute tolerance of $10^{-10}$ and a relative tolerance of $10^{-6}$.

\begin{figure}[t]
\centering
\includegraphics{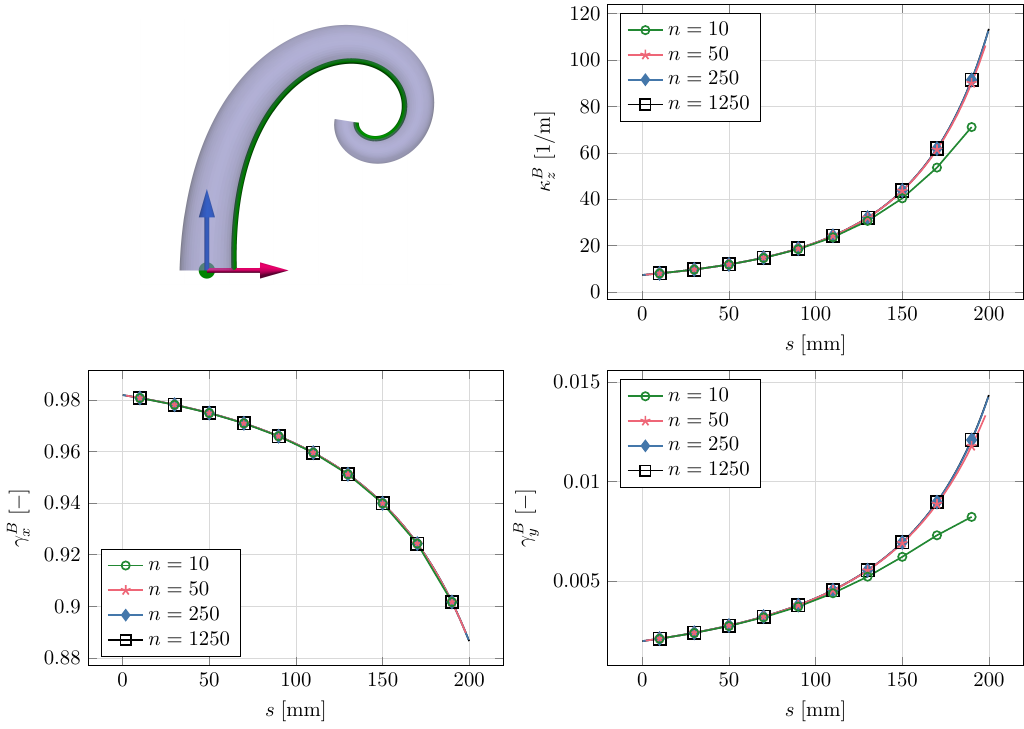}
\caption{Tendon-driven continuum manipulator at full actuation $\lambda_t=\SI{4}{\newton}$. Top-left: deformed configuration of the tapered rod ($r_\mathrm{tip}/r_0=0.4$) actuated by a single lateral tendon (green), producing a large spiral deformation. Top-right: bending curvature $\kappa_z^B$ along the rod centerline. Bottom-left and bottom-right: axial strain $\gamma_x^B$ and shear strain $\gamma_y^B$, respectively. Results are shown for four spatial discretizations ($n=10$, $50$, $250$, and $1250$). For the finest discretization, $\kappa_z^B$ increases from $\SI{7.3}{\meter\inv}$ to $\SI{113.4}{\meter\inv}$, $\gamma_x^B$ decreases from $0.982$ to $0.887$, and $\gamma_y^B$ increases from $0.002$ to $0.014$. The responses for $n=50$, $250$, and $1250$ nearly overlap, indicating rapid mesh convergence.}
\label{fig:tdcrobot}

\end{figure}

Figure~\ref{fig:tdcrobot} shows the deformed configuration together with the distributions of the discrete strain measures along the rod centerline for four different spatial discretizations ($n=10$, $50$, $250$, and $1250$). The lateral tendon offset generates a distributed bending moment that drives the rod into a large spiral-like deformation. The coarsest discretization
($n=10$) exhibits visible deviations in the high-curvature tip region, whereas the results for $n=50$, $250$, and $1250$ are nearly indistinguishable, demonstrating rapid mesh convergence of the proposed discrete rod formulation. For the finest discretization, the axial strain~$\gamma_x^B$ decreases monotonically from $0.982$ at the clamped end to $0.887$ at the free tip, reflecting a moderate compressive axial load induced by the tendon tension along the deformed centerline. The bending curvature~$\kappa_z^B$ increases monotonically from $\SI{7.3}{\meter\inv}$ at the clamped end to $\SI{113.4}{\meter\inv}$ at the free tip, and the shear strain~$\gamma_y^B$ similarly increases from $0.002$ at the clamped end to $0.014$ at the free tip. This pronounced variation is a direct consequence of the tapered cross section. Since the axial, bending, and shear stiffnesses $EA(s) \propto r(s)^2$, $EI(s) \propto r(s)^4$, $GA(s) \propto r(s)^2$  decrease rapidly toward the tip, the same tendon tension produces progressively larger compressive axial strain, curvature, and shear deformation in the tip region, consistent with the spiral shape visible in the figure. The discrete rod formulation handles this continuously varying stiffness naturally through the element-wise compliance law~\eqref{eq:compliance_sys}, with no special treatment required at any location along the rod.

\paragraph{Wilberforce pendulum}
To demonstrate the capability of our formulation for highly dynamic problems, we consider the Wilberforce pendulum~\cite{Harsch2021}, consisting of a helical spring clamped at its upper end with a steel cylinder attached at the lower free end. Four screws with adjustable nuts allow tuning of the cylinder's moment of inertia. The spring is made of spring steel EN~10270-1 with mass density $\rho = \SI{7850}{\kilogram/\meter^3}$, Young's modulus $E = \SI{206e9}{\pascal}$, and shear modulus $G = \SI{81.5e9}{\pascal}$, formed into a perfect helix with $n_c = 20$ coils, coil radius $R = \SI{16}{\milli\meter}$, wire diameter $d = \SI{1}{\milli\meter}$, and unloaded pitch $c = \SI{1}{\milli\meter}$. The reference configuration is constructed by placing the rod nodes on the exact helix such that the nodal orientations coincide with the Serret--Frenet frame~\cite{Harsch2021}. The spring is discretized with $n = 800$ elements (40 elements per coil). The pendulum bob (cylinder together with screws and nuts) is modeled as a rigid steel cylinder of radius $R_b = \SI{25}{\milli\meter}$, height $h = \SI{34}{\milli\meter}$ and mass $m_\mathrm{bob} = \pi \rho R_b^2 h $, rigidly attached to the lower end of the spring.

The static equilibrium under gravity is first computed via a Newton--Raphson method with 10 load increments, taking into account the distributed weight of the helical spring, the gravitational force on the bob, and an additional downward pulling force of magnitude $0.3\,m_\mathrm{bob}\,g$ with $g = \SI{9.81}{\meter/\second^2}$ applied to the bob to produce an initial vertical displacement. Starting from this statically deformed configuration, the pulling force is removed and the dynamic simulation is performed over $T = \SI{20}{\second}$, with the system released from rest. The bob then oscillates longitudinally while geometric coupling in the deformed helix induces a torsional oscillation. The dynamic equations of motion~\eqref{eq:kin_diff_eq_sys}--\eqref{eq:compliance_sys} are integrated using a Radau-type solver for DAEs~\cite{DAESolver}, with the corresponding code published in~\cite{GitRepoDAESolver}. The simulations are performed using the default solver settings, with an absolute tolerance of $10^{-6}$ and a relative tolerance of $10^{-3}$.

\begin{figure}[htbp]
\centering

% ================= LEFT IMAGE =================
\begin{minipage}[c]{0.2\textwidth}
    \includegraphics[width=\linewidth]{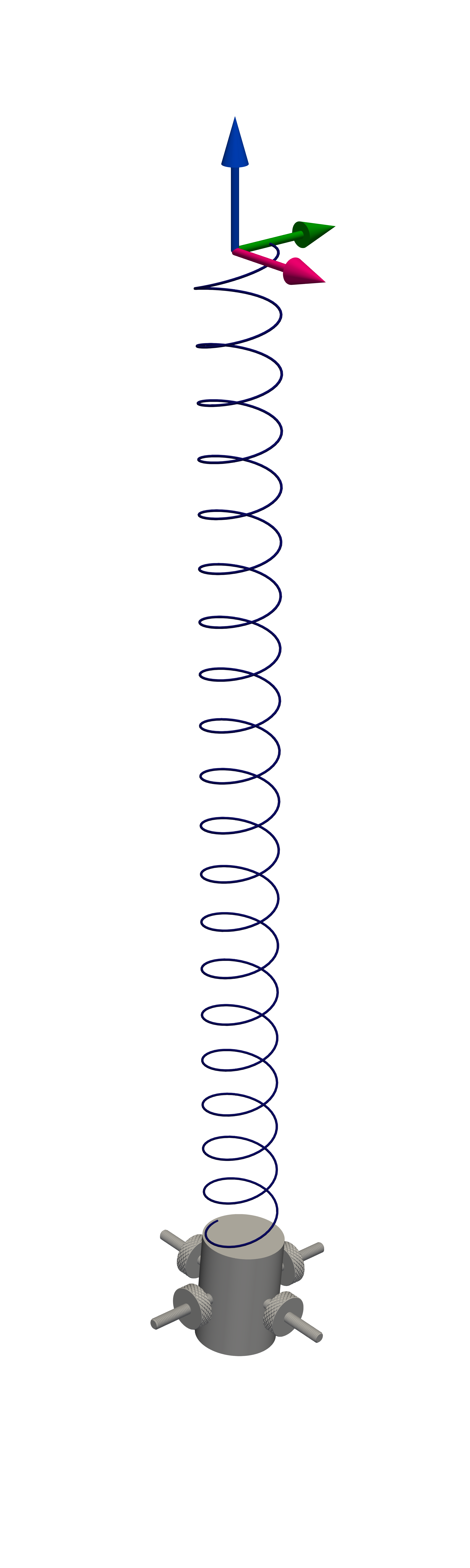}
\end{minipage}
% ================= RIGHT PLOTS =================
\begin{minipage}[c]{0.79\textwidth}
	\includegraphics{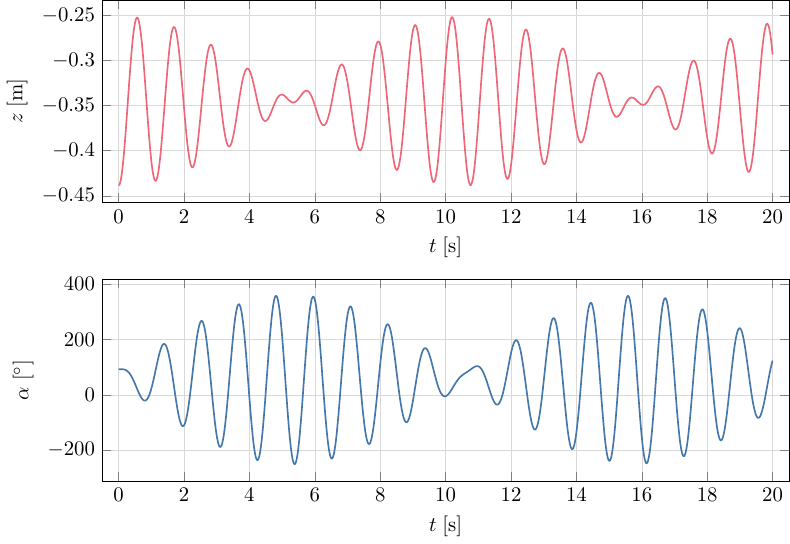}
%    \begin{tikzpicture}
%    \begin{groupplot}[
%        group style={
%            group size=1 by 2,
%            vertical sep=1.3cm
%        },
%        width=0.95\linewidth,
%        height=5cm,
%        grid=both,
%        grid style={line width=0.2pt, draw=gray!30},
%        tick align=inside,
%    ]
%
%    % ---------- Top plot ----------
%    \nextgroupplot[
%        ylabel={$z\;[\si{\meter}]$},
%        xlabel={$t\; [\si{\second}]$},
%		enlarge x limits={abs=0.5},
%		% enlarge y limits={abs=0.03},
%    ]
%
%    \addplot[
%        myred,
%        thick
%    ]
%    table[
%        col sep=comma,
%        x=t,
%        y=z
%    ] {figures/data_wilberforce2p0.csv};
%
%    % ---------- Bottom plot ----------
%    \nextgroupplot[
%        xlabel={$t\;[\si{\second}]$},
%        ylabel={$\alpha\;[\si{\degree}]$},
%		enlarge x limits={abs=0.5},
%		% enlarge y limits={abs=70},
%    ]
%
%    \addplot[
%        myblue,
%        thick
%    ]
%    table[
%        col sep=comma,
%        x=t,
%        y=alpha
%    ] {figures/data_wilberforce2p0.csv};
%
%    \end{groupplot}
%    \end{tikzpicture}
\end{minipage}
\caption{Wilberforce pendulum example. Left: deformed configuration of the helical spring with the attached cylindrical bob. Right: time history of the vertical displacement $z(t)$ of the bob center of mass (top) and the rotation angle $\alpha(t)$ about the spring axis (bottom) over $T=\SI{20}{\second}$.}
\label{fig:wilberforce}

\end{figure}

The time histories of the vertical position $z(t)$ and the rotation angle $\alpha(t)$ of the bob are shown in Figure~\ref{fig:wilberforce}. The vertical displacement oscillates between $\SI{-0.252}{\meter}$ and $\SI{-0.439}{\meter}$, while the torsional angle $\alpha$ reaches peak amplitudes of approximately $\pm 300^\circ$. The characteristic energy exchange between the two modes is clearly visible. As the vertical amplitude decreases, the torsional amplitude grows, and vice versa, consistent with the behavior reported in~\cite{Harsch2021}. With the cylinder's moment of inertia appropriately tuned, the two modes exhibit an almost perfect phase shift of $\pi/2$. This result demonstrates that the proposed discrete rod formulation correctly captures the geometrically induced coupling between longitudinal and torsional deformation modes in a dynamically loaded helical spring.

\section{Conclusion} \label{sec:conclusion}
We have presented a mixed discrete Cosserat rod formulation that models a slender elastic rod as a chain of rigid bodies (nodes) coupled by compliant elastic forces and moments acting between adjacent node pairs. Discrete dilatation-shear and torsion-curvature strain measures of the rod are computed directly from the relative positions and orientations of each node pair, and the constitutive behavior is described by element-wise compliance laws consistent with the Hellinger--Reissner variational principle. This discrete multibody representation emerges from the mixed Petrov--Galerkin Cosserat rod FEM of Herrmann et al.~\cite{Herrmann2025} at polynomial degree $p=1$, combined with the inertial virtual work framework of Harsch et al.~\cite{Harsch2023}, where the internal virtual work is evaluated at element midpoints by the midpoint rule, and the external and inertial virtual work at nodes by the trapezoidal rule. The resulting model exposes the same two-node kinematic topology as discrete rod models from the computer graphics community, while inheriting the robustness and the absence of locking from the underlying mixed FEM. This locality of the force and moment interactions, together with the locking-free behavior, makes the formulation a promising starting point for GPU-accelerated implementations of Cosserat rod FEMs.

The three numerical examples demonstrate the capability of the proposed formulation across a range of geometric complexities and loading conditions. The two-helical-segments benchmark confirms that the element-wise compliance law correctly handles piecewise-varying cross sections without any special treatment at cross-sectional transitions. The tendon-driven continuum manipulator shows that the formulation naturally handles continuously varying cross-sectional geometry and large spiral deformations while exhibiting rapid mesh convergence across different spatial discretizations. The monotonically increasing bending curvature toward the tip is fully consistent with the decreasing bending stiffness $EI(s) \propto r(s)^4$ of the tapered cross-sectional geometry. The Wilberforce pendulum example demonstrates the accuracy of the dynamic formulation, reproducing the characteristic energy exchange between longitudinal and torsional oscillation modes with an almost perfect phase shift of $\pi/2$.

Future work will address (i) the design of local relaxation solvers compatible with the mixed saddle-point structure of the present formulation, enabling GPU implementations that exploit the two-node connectivity, and (ii) the extension to frictional contact interaction between several rods and surrounding objects.

\begin{acknowledgement}
This research was funded by the Deutsche Forschungsgemeinschaft (DFG, German Research Foundation; grant number 405032572) as part of the DFG Priority Programme Soft Material Robotic Systems (SPP2100). We also thank Philipp L. Kinon (Karlsruhe Institute of Technology) for fruitful discussions.
\end{acknowledgement}

\def\bstname{pamm}
% Use this code if you wish to generate your bibliography with BibTeX;
% please replace first the string "demo" below with the name(s) of
% the BibTeX data base(s) you want to use.
% The resulting bibliography-output (the contents of the .bbl file)
% must be pasted into this file before submission.
% 
% \bibliographystyle{pamm}
% \bibliography{references.bib}
% 
% Replace the following example bibliography with your references
% before submission:
\providecommand{\WileyBibTextsc}{}
\let\textsc\WileyBibTextsc
\providecommand{\othercit}{}
\providecommand{\jr}[1]{#1}
\providecommand{\etal}{~et~al.}

\end{document}